\documentclass[12 pt]{article}%
\usepackage{amsmath, amsfonts, amsthm, color,latexsym}
\usepackage{amsmath, ulem}
\usepackage{amsfonts}
\usepackage{amssymb}
\usepackage{color, soul}
\usepackage{xcolor}
\usepackage[all]{xy}
\usepackage{graphicx}%
\usepackage{hyperref}
\providecommand{\U}[1]{\protect\rule{.1in}{.1in}}
\allowdisplaybreaks[4]
\newtheorem{theorem}{Theorem}[section]
\newtheorem{proposition}[theorem]{Proposition}
\newtheorem{corollary}[theorem]{Corollary}
\newtheorem{example}[theorem]{Example}

\newtheorem{remark}[theorem]{Remark}

\newtheorem{lemma}[theorem]{Lemma}
\newtheorem{final remark}[theorem]{Final Remark}

\allowdisplaybreaks[4]

\newcommand {\R}{\mathbb{R}}

\newcommand {\N} {\mathbb{N}}
\newcommand{\norma}[1]{\| #1 \|}
\newcommand{\conj}[2]{\left \{ {#1} \, : \, {#2} \right \}}

\begin{document}

\title{Arens extensions of disjointness preserving multilinear operators}
\author{Geraldo Botelho\thanks{Supported by FAPEMIG grants RED-00133-21 and APQ-02622-26.} , Luis Alberto Garcia\thanks{Supported by FAPEMIG grant APQ-02622-26.} ~and  Vinícius C. C. Miranda\thanks{Supported by FAPESP (2025/08630-0 and 2023/12916-1) and FAPEMIG (APQ-02622-26) \newline 2020 Mathematics Subject Classification: 46A40, 46B24, 46G25, 47H60.
  }}
\date{}
\maketitle

\begin{abstract}
First we give a counterexample showing that, unlike the second adjoint of a disjointness preserving linear operator, Arens extensions of disjointness preserving multilinear operators are not disjointness preserving in general. Then we prove that, if $E_1, \ldots, E_m$ are 
Riesz spaces with $F$ Archimedean and $A \colon E_1 \times \cdots \times E_m \to F$ is a regular disjointness preserving $m$-linear operator, then all Arens extensions of $A$ are disjointness preserving if: (i) $A$ has finite lattice rank; or (ii) $E_1, \ldots, E_m$ are also Archimedean and  the extensions are restricted to the product of the order continuous biduals; or (iii) the spaces are Banach lattices and either $E_1^*, \ldots, E_m^*$ have order continuous norms or $F^*$ has a Schauder basis consisting of disjointness preserving functionals. 
\end{abstract}

\section{Introduction}

The study of disjointness preserving linear operators between  Riesz spaces (vector lattices) or Banach lattices started with the seminal work of Abramovich, Veksler and Koldunov  \cite{abramovich} in 1979. Next, several authors investigated different aspects of such class of operators, see, e.g., \cite{arenson, kitover, boulabiar2008, kalauch, kuku, pagter}. Multilinear disjointness preserving operators were considered later by Kusraev and Tabuev \cite{tabuev} in 2004, and then studied by Bu, Buskes and Kusraev in \cite{buservey},  by Kusraev and Kusraeva in \cite{kus}, by Toumi in \cite{toumi}, and, more recently, by Kusraeva in \cite{kus2021, kus2023}. Following this trend, the purpose of this paper is to initiate the study of disjointness preserving multilinear operators whose Arens extensions are disjointness preserving as well.

To state the general problem, let us fix some notation. By $E^\sim$ we denote the order dual of the Riesz space $E$, hence $E^{\sim\sim}$ stands for its second order dual; and $J_E \colon E \to E^{\sim\sim}$ denotes the canonical Riesz homomorphism.   Let $E_1, \ldots, E_m,F$ be Riesz spaces, let $A \colon E_1 \times \cdots \times E_m \to F$ be a regular $m$-linear operator, and let $\rho$ be a permutation of the set $\{1, \ldots, m\}$. By
\begin{equation}\label{p9ms}AR_{m}^{\rho}(A)\colon E_{1}^{\sim\sim}\times\cdots\times E_{m}^{\sim\sim}\longrightarrow F^{\sim\sim}
\end{equation}
we mean the Arens extension of $A$ with respect to $\rho$, which happens to be a regular $m$-linear operator that extends $A$ in the sense that
\begin{equation}\label{equaextensao}
    AR_{m}^{\rho}(A)(J_{E_{1}}(x_{1}),\ldots,J_{E_{m}}(x_{m}))=J_{F}(A(x_1,\ldots,x_m))
\end{equation}
for all $x_1\in E_1,\ldots,x_m\in E_m$. Details can be found in Section 3. In particular, Arens extensions are the multilinear analogues of the second adjoint of a linear operator.

\medskip

{\bf Problem.} Is $AR_{m}^{\rho}(A)$ disjointness preserving whenever $A$ is?

\medskip

The corresponding problems for Riesz multimorphisms, for almost Dunford-Pettis multilinear operators and for order continuous multilinear operators were addressed by the first two authors, respectively, in \cite{geraldoluis, gl1, edinburgh}.

 In Section \ref{contra}, we will present a counterexample that solves negatively the problem above, as well as problem (ii) posed in \cite[p.\,2]{geraldoluis}. The counterexample gives a negative solution to the problem even in the more restricted environment of Banach lattices. The next obvious step is to find conditions under which the problem has a positive solution. We shall give contributions in the two following directions: 

(I) For Riesz subspaces $X_j$ of $E_j^{\sim\sim}$ containing $J_{E_j}(E_j)$, $j = 1, \ldots, m$, the restriction of $AR_{m}^{\rho}(A)$ to $X_1 \times \cdots \times X_m$ shall be denoted by $AR_{m}^{\rho}(A)\colon X_1 \times \cdots \times X_m \to F^{\sim\sim}$. It is easy to see that, if $A$ preserves disjointness, then so does
$$AR_{m}^{\rho}(A)\colon J_{E_1}(E_1)\times\cdots\times J_{E_m}(E_{m})\longrightarrow F^{\sim\sim}.$$
We shall prove positive results for some products $X_1 \times \cdots \times X_m$ larger than $J_{E_1}(E_1)\times\cdots\times J_{E_m}(E_{m})$ and smaller (in general) than $E_{1}^{\sim\sim}\times\cdots\times E_{m}^{\sim\sim}$, e.g., Theorem \ref{proporest}.

(II) Positive results for the whole product $E_{1}^{\sim\sim}\times\cdots\times E_{m}^{\sim\sim}$ shall be obtained under conditions on the operator $A$ (Theorem \ref{teo2-6}), or on the domain spaces $E_1, \ldots, E_m$ (Corollary \ref{p09j}), or on the target space $F$ (Corollary \ref{cor2}), or on the domain and target spaces (Corollary \ref{cor4}).

Denoting by $E^*$ the dual of a Banach lattice $E$, when the spaces above are Banach lattices, each $E_j^{\sim\sim}$ and $F^{\sim\sim}$ should be replaced with $E_j^{**}$ and $F^{**}$ and the canonical Riesz homomorphisms $J_{E_j}$ and $J_F$ becomes the canonical isometric embedding.

Since, as far as we know, this is the first time this problem is being addressed, some preparation is needed before tackling the main problem. It just so happens that the linear case of the problem has not been investigated yet. For further use, we prove in Section 2 that if $T \colon E \to F$ is an order bounded disjointness preserving linear operator between Riesz spaces with $F$ Archimedean, then its second adjoint $T^{\sim\sim} \colon E^{\sim\sim} \to F^{\sim\sim}$ preserves disjointness as well. For operators between Banach lattices no assumption on $F$ is required. Preparatory results on disjointness preserving multilinear operators are  proved in Section 3. It is worth mentioning that, in several results on operators between Riesz spaces, the domain spaces are not supposed to be Archimedean as they were in \cite{boulabiar2008, kus, kus2021, pagter}.

Section 4 is devoted to the aforementioned counterexample.

Let us summarize the main results proved in Section 5. With the notation fixed in (\ref{p9ms}), the following hold for a disjointness preserving operator $A$: (i) If the target Riesz space is Archimedean and $A$ has finite lattice rank, then $AR_{m}^{\rho}(A)$ is disjointness preserving. (ii) If all spaces are Archimedean, then the restriction of $AR_{m}^{\rho}(A)$ to the product $(E_{1}^{\sim})_{n}^{\sim}\times\cdots\times (E_{m}^{\sim})_{n}^{\sim}$ of the order continuous biduals is disjointness preserving. (iii) Supposing that the spaces are Banach lattices, $AR_{m}^{\rho}(A)$ is disjointness preserving on $E_1^{**} \times \cdots \times E_m^{**}$ if either $E_1^*, \ldots, E_m^*$ have order continuous norms; or $F^*$ has a Schauder basis consisting of disjointness preserving functionals; or $E$ is Arens regular, $E^*$ is separable and the lattice operations on $F^{**}$ are weak$^*$-sequentially continuous. 

We refer the reader to \cite{alip, meyer} for background on Riesz spaces and Banach lattices. Let us fix some terminology for multilinear maps between Riesz spaces.
Given Riesz spaces $E_1, \ldots, E_m,F$, an $m$-linear operator $A \colon E_1 \times \cdots \times E_m \to F$ is {\it positive}, written $A \geq 0$, if
$$A(x_1, \ldots, x_m) \geq 0 \mbox{~for all~} 0 \leq x_1 \in E_1, \ldots, 0 \leq x_m \in E_m. $$
An operator is {\it regular} if it can be written as the difference of two positive operators, and the space of all such operators is denoted by ${\cal L}_r(E_1, \ldots, E_m;F)$. For precise definitions, see \cite{buservey, loane}. Riesz spaces and Banach lattices are supposed to be real. When $F = \R$, we write  ${\cal L}_r(E_1, \ldots, E_m)$.

The positive cone of a Riesz space $E$ is denoted by $E^+$.
As usual, the order dual of a Riesz space $E$ is denoted by $E^\sim$ and the topological dual of a Banach lattice $E$ by $E^*$. The symbol $x \perp y$ means that the vectors $x$ and $y$ of a Riesz space are disjoint. The set of all permutations of $\{1, \ldots, m\}$ is denoted by  $S_m$.

\section{Disjointness preserving linear operators}

As mentioned in the Introduction, this section and the next one are a preparation for the results we prove in sections 4 and 5.

Recall that the (order) adjoint of an order bounded linear operator $T \colon E \to F$ between two Riesz spaces is the restriction of the adjoint of $T$ to $F^\sim$, denoted by $T^\sim\colon F^\sim \to E^\sim$. Of course, $T^{\sim\sim}: = (T^\sim)^\sim$.
As announced in the Introduction, we begin this section of preparatory results studying when the second adjoint $T^{\sim\sim} \colon E^{\sim \sim} \to F^{\sim \sim}$ of a disjointness preserving linear operator $T \colon E \to F$ between two Riesz spaces is disjointness preserving as well.
Let us recall that a linear operator $T \colon E \to F$ between two Riesz spaces is said to be disjointness preserving if $T(x) \perp T(y)$ in $F$ whenever $x \perp y$ in $E$. 

We begin our investigation considering order bounded linear operators between Riesz spaces. For a Riesz space $E$, by $J_E \colon E \to E^{\sim\sim}$ we denote the canonical Riesz homomorphism given by $J_E(x)(x^\sim) = x^\sim(x)$ for $x \in E$ and $x^\sim \in E^\sim$. By \cite[Theorem 1.69]{alip}, $J_{E}$ is injective whenever $E^{\sim}$ separates the points of $E$.


\begin{theorem} \label{teo1}
    Let $T \colon E \to F$ be an order bounded linear operator between two Riesz spaces. \\
    {\rm (1)} If $F$ is Archimedean and $T$ is disjointness preserving, then $T^{\sim\sim} \colon  E^{\sim \sim} \to F^{\sim \sim}$ is disjointness preserving. \\
    {\rm (2)} If 
    $J_{E}$ is a bijection and $T$ is disjointness preserving, then $T^{\sim\sim} \colon E^{\sim \sim} \to F^{\sim \sim}$ is disjointness preserving. \\
    {\rm (3)} If $F^{\sim}$ separates the points of $F$ and  $T^{\sim\sim} \colon E^{\sim \sim} \to F^{\sim \sim}$ is disjointness preserving, then $T$ is disjointness preserving.
\end{theorem}

\begin{proof}
    (1)  By \cite[Theorem 2.40]{alip} the modulus of $T$ exists and 
     $|T|(|x|)=|T(|x|)|=|T(x)|$ for every $x \in E$. Applying \cite[Theorem 1.7]{alip} we get
    \begin{align*}
        ||T|(x)|&=||T|(x^{+})-|T|(x^{-})|=||T(x^{+})|-|T(x^{-})||\\
        &=|T(x^{+})+T(x^{-})|\wedge |T(x^{+})-T(x^{-})|=|T(|x|)|\wedge |T(x)|=|T|(|x|),
    \end{align*}
hence $ ||T|(x)|=|T|(|x|)$ for every $x\in E$. This yields that $|T|$  is Riesz homomorphism, therefore the operator $|T|^{\sim\sim}$ is also a Riesz homomorphism by \cite[Theorems 2.19 and 2.20]{alip}. Given  $x^{\sim\sim}, y^{\sim\sim} \in E^{\sim\sim}$ such that $x^{\sim\sim}\perp y^{\sim\sim}$, we have
    \begin{align*}
        0&\leq |T^{\sim\sim}(x^{\sim\sim})|\wedge |T^{\sim\sim}(y^{\sim\sim})|\leq |T^{\sim\sim}|(|x^{\sim\sim}|)\wedge |T^{\sim\sim}|(|y^{\sim\sim}|)\\
        &\leq |T|^{\sim\sim}(|x^{\sim\sim}|)\wedge |T|^{\sim\sim}(|y^{\sim\sim}|)=|T|^{\sim\sim}(|x^{\sim\sim}|\wedge |y^{\sim\sim}|)=|T|^{\sim\sim}(0)=0,
    \end{align*}
which implies that $T^{\sim\sim}(x^{\sim\sim})\perp T^{\sim\sim}(y^{\sim\sim})$, proving that $T^{\sim\sim}$ is disjointness preserving.

(2) Since $T^{\sim\sim}(J_E(x)) = J_F(T(x))$ for every $x \in E$, we get $T^{\sim\sim}(x) = J_F(T(x))$ for every $x \in J_E(E) = E^{\sim\sim}$, which gives that $T^{\sim\sim}$ is disjointness preserving because both $T$ and $J_F$ are disjointness preserving.

(3) Let $x$ and $y$ be disjoint vectors in $E$. As $J_E \colon E \to E^{\sim \sim}$ is a Riesz homomorphism, we have $J_E(x) \perp J_F(y)$, from which it follows that
$$  J_F(T(x)) = T^{\sim\sim}(J_E(x))\perp T^{\sim\sim}(J_E(y)) = J_F(T(y)) \mbox{ in }F^{\sim\sim}$$
because $T^{\sim\sim}$ preserves disjointness. Using now that $(J_F)^{-1}\colon J_F(F) \to F$ is a Riesz homomorphism and that $F^\sim$ separates the points of $F$, we get $T(x) \perp T(y)$.  
\end{proof}

If $F$ is a Banach lattice, then $F$ is Archimedean and its topological dual $F^* = F^\sim$ separates the points of $F$. Therefore, the following is a straightforward consequence of Theorem \ref{teo1}.

\begin{corollary} \label{cor1}
    Let $E$ be a Riesz space and let $F$ be a Banach lattice. Then, an order bounded linear operator $T\colon E \to F$ is disjointness preserving if and only if its second adjoint $T^{\sim\sim} \colon E^{\sim \sim} \to F^{**}$ is disjointness preserving.
\end{corollary}


\begin{remark}\rm
The situation is quite simpler for positive operators. The following are equivalent for a positive linear operator $T \colon E \to F$ between Riesz spaces, where  $F^{\sim}$ separates the points of $F$.\\
\indent{\rm (a)} $T$ is disjointness preserving.\\
\indent{\rm (b)} $T$ is a Riesz homomorphism.\\
\indent{\rm (c)} $T^{\sim\sim} \colon E^{\sim \sim} \to F^{\sim \sim}$ is disjointness preserving.\\
\indent{\rm (d)} $T^{\sim\sim} \colon E^{\sim \sim} \to F^{\sim \sim}$  is a Riesz homomorphism. \\
 Indeed, (a)$\Leftrightarrow$(b) is well known, and it is easy to be checked directly; so (c)$\Leftrightarrow$(d) also holds because $T^{\sim\sim} \geq 0$. (a)$\Rightarrow$(d) follows from a combination of \cite[Theorem 2.19]{alip}, which imposes no conditions on $E$ and $F$, with the implication of \cite[Theorem 2.20]{alip} that does not require any assumption on the spaces. For (d)$\Rightarrow$(a), use that $E$ and $F$ are (Riesz isomorphic to) sublattices of $E^{\sim\sim}$ and $F^{\sim\sim}$, or, alternatively, that $T = (J_F)^{-1} \circ T^{\sim\sim} \circ J_E$, where $(J_F)^{-1} \colon J_F(F) \to F$.
\end{remark}

\section{Disjointness preserving multilinear operators}

As in the linear case, a multilinear operator $A\colon E_1 \times \cdots \times E_m \to F$ between Riesz spaces is said to be {\it order bounded} if $A$ sends order bounded subsets of $E_1 \times \cdots \times E_m$, with respect to the coordinatewise ordering, to order bounded subsets of $F$. It is easy to check that every regular $m$-linear operator is order bounded.

The symbol $\stackrel{\hat{j}}{\ldots}$ means that the $j$-th coordinate has been omitted. Given Riesz spaces $E_1, \dots, E_m,$ and  $ F$, an order bounded $m$-linear operator $A\colon E_1 \times \cdots \times E_m \to F$ is said to {\it preserve disjointness} if $A$ preserves disjointness separately, that is, for all $j\in\{1,\ldots,m\}$ and $x_{i}\in E_{i},\, i\neq j$, the linear operator
   $$A_{x_1,\stackrel{\hat{j}}{\ldots} ,x_m}\colon E_j \to F~,~ 
   A_{x_1,\stackrel{\hat{j}}{\ldots} ,x_m}(x_j) = A(x_{1}\ldots,x_{m}),$$
    preserves disjointness.
We shall apply, multiple times, a result that can be found in \cite[Theorem 3]{kus}, where the reader is referred to \cite{tabuev} for a proof. Since \cite{tabuev} is a difficult-to-read reference (it is in Russian), next we provide a proof for the benefit of the reader.


Recall that an $m$-linear operator $A\colon E_1 \times \cdots \times E_m \to F$ between Riesz spaces is a {\it Riesz multimorphism} (or an {\it $m$-morphism}) if $|A(x_1, \ldots, x_m)| = A(|x_1|, \ldots, |x_m|)$ for all $x_1 \in E_1,\ldots, x_m \in E_m$. Equivalently, $A$ is a Riesz multimorphism if and only if it is a Riesz homomorphism separately, that is, for every $j$, the linear map obtained by fixing positive vectors in all variables different from the $j$-th one is a Riesz homomorphism. 

\begin{remark}\label{4tui}\rm
 Clearly, 
Riesz multimorphisms preserve disjointness. However, the converse is not true in general: the map 
$$A\colon \ell_2\times\ell_2\longrightarrow\mathbb{R}~,~ A\Large{(}(x_n)_n,(y_n)_n \Large{)}=-x_1y_1,$$ 
is a disjointness preserving regular bilinear form that fails to be a 2-morphism. As in the linear case, a positive multilinear operator preserves disjointness if and only if it is a Riesz multimorphism.
\end{remark}


\begin{proposition}\label{prop2} Let $E_{1},\ldots,E_{m}, F$ be Riesz spaces with $F$ Archimedean and let $A \colon E_{1}\times\cdots \times E_{m} \to F$ be an order bounded $m$-linear operator. If $A$ is disjointness preserving, then $|A|$, $A^{+}$ and $A^{-}$ exist and they are Riesz multimorphisms.  Moreover,
    $$|A(|x_{1}|,\ldots,|x_{m}|)|=||A|(x_{1},\ldots,x_{m})|=|A(x_{1}\ldots,x_{m})| = |A|(|x_{1}|,\ldots,|x_{m}|)$$
 for all $x_{1}\in E_{1},\ldots,x_{m}\in E_{m}$, and
$$A^{+}(x_{1},\ldots,x_{m})=(A(x_{1},\ldots,x_{m}))^{+} \mbox{ and~} A^{-}(x_{1},\ldots,x_{m})=(A(x_{1},\ldots,x_{m}))^{-}$$
for all $x_{1}\in E_{1}^{+},\ldots,x_{m}\in E_{m}^{+}$.
\end{proposition}

\begin{proof} For simplicity, we prove the bilinear case $m=2$. Define $B_{1}\colon E_1^+ \times E_2^+ \to F^+$ by $B_{1}(x, y) = |A(x, y)|$. To check that $B_{1}$ is additive in the first variable, let $x_1, x_2 \in E_1^+$ and $y \in E_{2}^+$ be given. We have
    \begin{align*}
        B_{1}(x_1, y) + B_{1}(x_2, y) & = |A(x_1, y)| + |A(x_2, y)| = |A_y(x_1)| + |A_y(x_2)|,
    \end{align*}
    where $A_y \colon E_1 \to F$ is defined by $A_y(x) = A(x,y)$ for every $x \in E_1$. Since the linear order bounded operator $A_y$ is  disjointness preserving by assumption, we have from \cite[Theorem 2.40]{alip} that its modulus exists and $|A_y(x)| = |A_y|(x)$ for every $x \in E_1^+$. Therefore,
    \begin{align*}
        B_{1}(x_1, y) + B_{1}(x_2, y) & = |A_y|(x_1) + |A_y|(x_2) = |A_y|(x_1 + x_2) \\
        & = |A_y(x_1 + x_2)| = |A(x_1 + x_2, y)| = B_{1}(x_1 + x_2, y).
    \end{align*}
    This proves that $B_{1}$ is additive in the first variable. The additivity in the second variable follows analogously. From \cite[Theorem 2.3]{loane}, the map $\widetilde{B_{1}}\colon E_1 \times E_2 \to F$ given by
    $$ \widetilde{B_{1}}(x, y) = B_{1}(x^+, y^+) - B_{1}(x^+, y^-) - B_{1}(x^-, y^+) + B_{1}(x^-, y^-)  $$
is a positive bilinear extension of $B_{1}$. Let us check that $\widetilde{B_{1}} = A \vee (-A)=|A|$. Indeed, clearly we have $\widetilde{B_{1}} \geq A$ and $\widetilde{B_{1}} \geq -A$. Now, if $B \geq A$ and $B \geq -A$, we have that $B(x,y) \geq \pm A(x,y)$ for all $x \in E_1^{+}$ and $y\in E_2^{+}$, and so $B(x,y) \geq |A(x,y)| = \widetilde{B_{1}}(x,y)$, proving that  $\widetilde{B_{1}}=A\vee (-A) = |A|$ exists. Given $x\in E_{1}^{+}$ and $y\in E_{2}^{+}$, we have $$|A|(x,y)=\widetilde{B_{1}}(x,y)=B_{1}(x,y)=|A(x,y)|.$$

Now, let $x\in E_1$ and $y \in E_2$ be given. From the sentence above, it holds
   $|A|(|x|,|y|)=|A(|x|,|y|)|$. Since $A$ is disjointness preserving and $y^{+}\perp y^{-}$,
\begin{align*}
B_{1}(x^+, y^+) \wedge B_{1}(x^+, y^-)=|A(x^+, y^+)| \wedge |A(x^+, y^-)|=0.
\end{align*}
It follows that $B_{1}(x^+, y^+) \wedge B_{1}(x^+, y^-)=0$. In a similar way, $B_{1}(x^-, y^+) \wedge B_{1}(x^-, y^-)=0$. Moreover,
\begin{align*}
|B_{1}(x^+, y^+) - B_{1}(x^+, y^-)|&\wedge |B_{1}(x^-, y^+) - B_{1}(x^-, y^-)|\\
&=(B_{1}(x^+, y^+) + B_{1}(x^+, y^-))\wedge (B_{1}(x^-, y^+) + B_{1}(x^-, y^-))\\
&=B_{1}(x^+, |y|)\wedge B_{1}(x^-, |y|)= |A(x^+, |y|)|\wedge |A(x^-, |y|)|=0.
\end{align*}
Thus, for all $x\in E_{1}$ and $y\in E_{2}$,
\begin{align*}
||A|(x,y)|&=|\widetilde{B_{1}}(x, y)| = |B_{1}(x^+, y^+) - B_{1}(x^+, y^-) - (B_{1}(x^-, y^+) - B_{1}(x^-, y^-))|\\
&=|B_{1}(x^+, y^+) - B_{1}(x^+, y^-)| + |(B_{1}(x^-, y^+) - B_{1}(x^-, y^-))|\\
&=B_{1}(x^+, y^+)+ B_{1}(x^+, y^-)+ B_{1}(x^-, y^+)+ B_{1}(x^-, y^-)\\
&=B_{1}(|x|,|y|)=\widetilde{B_{1}}(|x|, |y|)=|A|(|x|,|y|).
\end{align*}
This proves that $|A|$ is a Riesz bimorphism.

    Similarly, defining $B_{2}(x,y) := (A(x,y))^+$  and $B_{3}(x,y) := (A(x,y))^{-}$ for all $x \in E_1^+$ and $y\in E_2^+$, we get that these maps yield the existence of $A^+$ and $(A)^{-}$, that both of them are Riesz bimorphisms, and that $A^{+}(x,y)=(A(x,y))^{+}$ and $A^{-}(x,y)=(A(x,y))^{-}$ for all $x\in E_{1}^{+}$ and $y\in E_{2}^{+}$, from which the other desired equalities follow.
\end{proof}

%
%




%

The following is an immediate consequence of Proposition \ref{prop2}.
\begin{corollary}\label{corol}
Let $E_{1},\ldots,E_{m}, F$ be Riesz spaces with $F$ Archimedean and let $A \colon E_{1}\times\cdots \times E_{m} \to F$ be an order bounded $m$-linear operator. Then, $A$ is disjointness preserving if, and only if, for every $j\in\{1,\ldots,m\}$ and fixed positive vectors $x_{i}\in E_{i},\, i\neq j$, the linear operator $
    x_j \in  E_{j}\longrightarrow A(x_{1}\ldots,x_{m}) \in F$
    preserves disjointness.
\end{corollary}

\section{The counterexample}\label{contra}

By Corollary \ref{cor1}, the second adjoint of a disjointness preserving regular linear operator between Banach lattices is disjointness preserving too. In this section we show that this is no longer for Arens extensions of multilinear operators.


We shall use the notation from \cite{gl1} to introduce Arens extensions of multilinear operators. Let $E_1, \ldots, E_m,F$ be Riesz spaces and let $m\in\mathbb{N}$ and $\rho\in S_{m}$ be given. For  $i=1,\ldots,m-1$, and $x_{\rho(i)}^{\sim\sim}\in E_{\rho(i)}^{\sim\sim}$, consider  the linear operator
\begin{equation*}
\overline{x_{\rho(i)}^{\sim\sim}}\colon \mathcal{L}_{r}(E_{\rho(i)},\ldots,E_{\rho(m)})\longrightarrow \mathcal{L}_{r}(E_{\rho(i+1)},\ldots,E_{\rho(m)})~,~ \overline{x_{\rho(i)}^{\sim\sim}}(B)=x_{\rho(i)}^{\sim\sim}\circ B^{i},
\end{equation*}
where $B^{i}\colon E_{\rho(i+1)}\times\cdots\times E_{\rho(m)}\longrightarrow E_{\rho(i)}^{\sim}$ is given by
$$B^{i}(x_{i+1},\ldots,x_{m})(x_{i})=B(x_{i},x_{i+1},\ldots,x_{m}) \text{ for all } x_{j}\in E_{\rho(j)}, j=i,\ldots,m.$$
In the case $i=m$, we define $\overline{x_{\rho(m)}^{\sim\sim}}\colon E_{\rho(m)}^{\sim}\longrightarrow \mathbb{R}$ by $ \overline{x_{\rho(m)}^{\sim\sim}}=x_{\rho(m)}^{\sim\sim}$.

Given an $m$-linear form $C\colon E_{1}\times\cdots\times E_{m}\to \mathbb{R}$, to each permutation $\rho \in S_m$ we consider the $m$-linear form
$$C_{\rho}\colon E_{\rho(1)}\times\cdots\times E_{\rho(m)}\longrightarrow \mathbb{R}~,~ C_{\rho}(x_{1},\ldots,x_{m})=C(x_{\rho^{-1}(1)},\ldots,x_{\rho^{-1}(m)}).$$
The Arens extension of a regular $m$-linear operator
 $A\colon E_{1}\times\cdots\times E_{m}\longrightarrow F$ with respect to the permutation $\rho$ is defined as the $m$-linear operator $AR_{m}^{\rho}(A)\colon E_{1}^{\sim\sim}\times\cdots\times E_{m}^{\sim\sim}\longrightarrow F^{\sim\sim}$ given by
$$AR_{m}^{\rho}(A)(x_{1}^{\sim\sim},\ldots,x_{m}^{\sim\sim})(y^\sim)=\big(\overline{x_{\rho(m)}^{\sim\sim}}\circ\cdots\circ \overline{x_{\rho(1)}^{\sim\sim}}\big)((y^{\sim}\circ A)_{\rho})$$
for all $x_{1}^{\sim\sim}\in E_{1}^{\sim\sim},\ldots,x_{m}^{\sim\sim}\in E_{m}^{\sim\sim}$ and $y^{\sim}\in F^{\sim}$.
By \cite[Theorem 2.2(b)]{geraldoluis}, $AR_{m}^{\rho}(A)$ is an extension of $A$ in the sense that (\ref{equaextensao}) holds. 
Although the notation is a bit different, these are exactly the Arens extensions studied in \cite{geraldoluis} and in many papers on Arens (or Aron-Berner) extensions. The extension $A^{(m+1)*}$ from \cite{ryan1, Buskes} is recovered by considering the permutation
\begin{equation}\label{i8m3}\theta(m) = 1, \,  \theta(m-1) = 2, \ldots,\theta(2) = m-1, \, \theta(1) = m,
\end{equation}
 that is, $AR_m^\theta(A) = A^{*(m+1)}$.


Recall that the {\it ultrapower} of a Banach space $X$ with respect to the ultrafilter $\mathcal{U}$ on $\N$ is the quotient space $X_{\mathcal{U}} = \ell_\infty (X) / \mathcal{N}_{\mathcal{U}}$, where
$ \displaystyle \mathcal{N}_{\mathcal{U}} = \conj{(x_n)_n \in \ell_\infty}{\lim_{\mathcal{U}} \norma{x_n} = 0}. $
The class of $(x_n)_n$ in $X_{\mathcal{U}}$ is denoted by $[(x_n)_n]_{\mathcal{U}}$. We shall consider the  canonical isometric embedding
$$\iota_X \colon X \to X_{\mathcal{U}}~,~ \iota_X(x) = [(x, x, \dots)]_{\mathcal{U}},$$
and the quotient map
$$ Q_X \colon \ell_\infty(X) \to X_{\mathcal{U}}~,~Q_X((x_n)_n) = [(x_n)_n]_{\mathcal{U}}. $$
Moreover, there is also a canonical embedding
from the ultrapower of $X^*$ into the dual of the ultrapower of $X$:
$$j_X \colon (X^*)_{\mathcal{U}} \to (X_{\mathcal{U}})^*~,~
 \displaystyle j_X([(x_n^*)_n]_{\mathcal{U}})([(x_n)_n]_{\mathcal{U}}) = \lim_{\mathcal{U}}  x_n^*(x_n).$$
If $E$ is a Banach lattice, $E_{\mathcal{U}}$ is a Banach lattice with the order given pointwisely. We refer the reader to \cite{aviles, heinrich} for the theory of ultrapower of Banach spaces/lattices.

Let $E = \ell_\infty(\ell_1),$ $F = \ell_1(\ell_\infty)$ and $G = \ell_1(\ell_1)$. Consider the positive bilinear map $A: E \times F \to G$ defined by
$$ A(x,y) = (x_n \!\cdot  y_n)_{n=1}^\infty, \quad x = (x_n)_n \in E, \, y = (y_n)_n \in F, $$
where $x_n \cdot y_n$ denotes the pointwise product between $x_n \in \ell_1$ and $y_n \in \ell_\infty$.
Since $|A(x,y)| = A(|x|, |y|)$, $A$ is a Riesz bimorphism, hence $A$ is  disjointness preserving. We will split the proof that $A^{***} $ is not disjoitness preserving into two claims. To do this, we need to fix a free ultrafilter $\mathcal{U}$ on $\N$. The existence of free ultrafilters on $\N$ is consistent with ZFC and they exist by Zorn's Lemma. 

\medskip

\noindent {\bf Claim 1.} For every $x^{**} \in E^{**}$, $A^{***}(x^{**}, \eta) = R^*(x^{**})$, where
$\eta: \ell_\infty(\ell_\infty^*) \to \R$ is the positive linear functional given by
$$ \eta((\nu_n)_n) = \lim_{\mathcal{U}} \nu_n(\boldsymbol{1}), \quad (\nu_n)_n \in \ell_\infty(\ell_\infty^*), $$
where $\boldsymbol{1} = (1, 1, \dots) \in \ell_\infty, $
and $R : G^* \to E^*$ is the positive linear operator given by $$ R(x^{\ast})(x) = \lim_{\mathcal{U}} x_{n}^{\ast}(x_n), \quad x^{\ast} = (x^{\ast}_n)_n \in G^*, \, x = (x_n)_n \in E.$$
\noindent{\it Proof of Claim 1.} The easy verification that $\eta$ and $R$ are positive linear maps is omitted. 
To prove that $A^{***}(x^{**}, \eta) = R^*(x^{**})$ holds for every $x^{**} \in E^{**}$, it is enough to check that $\eta(y^{\ast}\circ A(x,\cdot))=R(y^{\ast})(x)$ holds for every $y^{\ast}\in G^{\ast}$ and every $x\in E$. To see this, consider
the Riesz isometric isomorphism
$$\psi\colon \ell_{\infty}(\ell_{1}^{\ast})\longrightarrow \ell_{1}(\ell_{1})^{\ast}~,~\psi((z_n^{\ast})_n)((z_n)_n)=\displaystyle\sum_{n=1}^{\infty} z_{n}^{\ast}(z_n),$$
$y^{\ast}\in G^{\ast}=\ell_{1}(\ell_{1})^{\ast}$, and $x=(x_n)_n\in E=\ell_{\infty}(\ell_{1})$. Then there exists $(z_{n}^{\ast})_n\in \ell_{\infty}(\ell_{1}^{\ast})$ such that $\psi((z_{n}^{\ast})_n)=y^{\ast}$.
For every $y = (y_n)_n \in F=\ell_{1}(\ell_{\infty})$, we have
\begin{align*}
    (y^* \circ A(x, \cdot))(y) &  = y^*(A(x,y)) = \psi((z_n^*)_n)(x_ny_n)_n = \sum_{n=1}^\infty z_n^*(x_ny_n).
\end{align*}
For each $n \in \N$, consider the positive linear operator $i_{n}\colon \ell_{\infty}\longrightarrow\ell_{1}(\ell_{\infty})$ given by $ i_{n}(y)=(0,\ldots,0,y,0\ldots)$. Consider also the Riesz isometric isomorphism
$$\phi\colon \ell_{\infty}(\ell_{\infty}^{\ast})\longrightarrow \ell_{1}(\ell_{\infty})^{\ast}~,~\psi((w_n^{\ast})_n)((w_n)_n)=\displaystyle\sum_{n=1}^{\infty} w_{n}^{\ast}(w_n).$$
Since $y^{\ast}\circ A(x,\cdot)\in \ell_{1}(\ell_{\infty})^{\ast}$, there exists $(w_{n}^{\ast})_{n}\in \ell_{\infty}(\ell_{\infty}^{\ast})$ such that $\phi((w_{n}^{\ast})_{n})=y^{\ast}\circ A(x,\cdot)$. Thus, for every $n\in\mathbb{N}$, $w_{n}^{\ast}=y^{\ast}\circ A(x,\cdot)\circ i_{n}$. Therefore,
\begin{align*}
    \eta(y^* \circ A(x, \cdot))&=\eta(\phi((w_{n}^{\ast})_{n}))=\lim_{\mathcal{U}} w_{n}^{\ast}(\textbf{1}) = \lim_{\mathcal{U}} (y^* \circ A(x, \cdot)) (i_n(\boldsymbol{1}))= \lim_{\mathcal{U}} y^*(A(x,i_n(\boldsymbol{1})))\\
    &= \lim_{\mathcal{U}} \psi((z_n^*)_n) (0, \dots, 0, x_n, 0, \dots) = \lim_{\mathcal{U}} z_n^*(x_n) = R(y^*)(x).
\end{align*}
Now, the definition of $A^{***}$
 yields that $A^{***}(x^{**}, \eta) = R^{*}(x^{**})$ for every $x^{**} \in E^{**}$, as claimed. ~$\square$

\medskip

\noindent {\bf Claim 2.} $R^*$ is not a Riesz homomorphism.

\medskip

\noindent{\it Proof of Claim 2.} First notice that
$ R = Q_{\ell_1}^* \circ j_{\ell_1} \circ Q_{\ell_\infty}$. Setting $e := j_{\ell_1}([\boldsymbol{1}]_{\mathcal{U}})$, we have that $e(y) = \norma{y}$ for every $y \in (\ell_1)_{\mathcal{U}}^+.$ Since $\mathcal{U}$ is a free ultrafilter on $\N$, $\mathcal{U}$ is completely incomplete \cite[p. 589]{aviles}. Thus, as $\ell_1$ is not super-reflexive, $j_{\ell_1}$ is not surjective \cite[Corollary 7.2]{heinrich}, and so there exists $0 < x^{\ast} \in ((\ell_1)_{\mathcal{U}})^*$ such that $ x^{\ast}\notin j_{\ell_1}((\ell_1^*)_{\mathcal{U}})$. So, for each $x \in (\ell_1)_{\mathcal{U}}^+$,
$$0 \leq x^{\ast}(x) \leq \norma{x^{\ast}} \norma{x} = \norma{x^{\ast}} e(x).$$
Hence, $\varphi := \dfrac{x^{\ast}}{\norma{x^{\ast}}} \in ((\ell_1)_{\mathcal{U}})^*$, $\varphi \notin j_{\ell_1}((\ell_1^*)_{\mathcal{U}})$ and $0 \leq \varphi \leq e$. As $Q_{\ell_1}^*$ is positive,
$$ 0 \leq Q_{\ell_1}^*(\varphi) \leq Q_{\ell_1}^*(e) = R(\boldsymbol{1}, \boldsymbol{1}, \dots), $$
and since $Q_{\ell_\infty}$ is surjective, $R(G^*) = Q_{\ell_1}^*(j_{\ell_1}((\ell_1^*)_{\mathcal{U}}))$,
which is a closed subspace of $E^*$ as both $j_{\ell_1}$ and $Q_{\ell_\infty}^*$ are isometries. Let us see that $Q_{\ell_1}^*(\varphi) \notin R(G^*)$. Indeed, if $Q_{\ell_1}^*(\varphi) = Q_{\ell_1}^*(j_{\ell_1}(\xi))$ for some $\xi \in (\ell_1^*)_{\mathcal{U}}$, the injectivevity of $Q_{\ell_1}^*$ would give $\varphi = j_{\ell_1}(\xi)$, contrary to the choice of $\varphi$. So, by the Hahn-Banach separation theorem, there exists $z^{**} \in E^{**}$ such that
$z^{**}|_{R(G^*)} = 0$ and $z^{**}(Q_{\ell_1}^*(\varphi)) > 0$. For every $z^{\ast} \in G^*$, $R^*(z^{**})(z^{\ast}) = z^{**}(R(z^{\ast})) = 0, $ hence
$$|R^*(z^{**})|(y^{\ast})=\sup\{|R^*(z^{**})(w^{\ast})|: |w^{\ast}|\leq y^{\ast}\}=\sup\{|z^{**}(R(w^{\ast}))|: |w^{\ast}|\leq y^{\ast}\}=0$$
for every $0\leq y^{\ast} \in G^*$. On the other hand,
 $$R^*(|z^{**}|)(\boldsymbol{1}, \boldsymbol{1}, \dots) = |z^{**}|(R(\boldsymbol{1}, \boldsymbol{1}, \dots)) \geq |z^{**}|(Q_{\ell_1^*} (\varphi)) \geq z^{**}(Q_{\ell_1^*} (\varphi)) > 0,$$
from which it follow that $R^*$ is not Riesz homomorphism.~$\square$

\medskip


Finally, applying Claim 2, there exists $z^{**}\in E^{**}$ such that $|R^*(z^{**})|\neq R^*(|z^{**}|)$.  By Claim 1, we have $|A^{***}(z^{**}, \eta)| \neq A^{***}(|z^{**}|, \eta)$, proving that $A^{\ast\ast\ast}$ is not a Riesz bimorphism, hence $A^{\ast\ast\ast}$ is not disjointness preserving because it is a positive bilinear operator (cf. Remark \ref{4tui}).



\section{Disjointness preserving Arens extensions}

Let $A\colon E_1\times\cdots\times E_m\longrightarrow F$ be a disjointness preserving regular $m$-linear operator between Riesz spaces, where each $E_{i}^\sim$ separates the points of $E_{i}$. 
For every $\rho \in S_m$, the fact that 
 $J_{F}$ is a Riesz homomorphism gives immediately that 
 $$AR_{m}^{\rho}(A)\colon J_{E_1}(E_1)\times\cdots\times J_{E_m}(E_{m})\longrightarrow F^{\sim\sim}$$
is disjointness preserving.

For a Riesz space $E$, by $(E^{\sim})_{n}^{\sim}$  we denote the subspace of $E^{\sim\sim}$ consisting of all order continuous functionals defined on $E^{\sim}$. It is worth noting that if $E^{\sim}$ separates the points of $E$, then $J_{E}(E)\subseteq (E^{\sim})_{n}^{\sim}\subseteq E^{\sim\sim}$ (see \cite[Theorem 1.69 and p. 61]{alip}). The first result of this section shows that Arens extensions preserve disjointness on the product of such restrictions.

\begin{theorem}\label{proporest}
    Let 
    $A\colon E_1\times\cdots\times E_m\longrightarrow F$ be a regular $m$-linear operator between Archimedean Riesz spaces. If $A$ is disjointness preserving, then all Arens extensions of $A$ are disjointness preserving on  $(E_{1}^{\sim})_{n}^{\sim}\times\cdots\times (E_{m}^{\sim})_{n}^{\sim}$. In other words, for every $\rho \in S_m$,
    $$AR_m^{\rho}(A) \colon (E_{1}^{\sim})_{n}^{\sim}\times\cdots\times (E_{m}^{\sim})_{n}^{\sim} \to F^{\sim\sim} $$
    is disjointness preserving.
\end{theorem}

\begin{proof} Let $\rho \in S_m$ be fixed. Given $i \in \{1,\ldots, m\}$, let $x_{i}^{\sim\sim}, z_{i}^{\sim\sim}\in (E_{i}^{\sim})_{n}^{\sim}$ be such that $x_{i}^{\sim\sim}\perp z_{i}^{\sim\sim}$. By Corollary \ref{corol} we can fix positive functionals $x_{j}^{\sim\sim}\in (E_{j}^{\sim})_{n}^{\sim}, j=1,\ldots,m$ with $j\neq i$. By Proposition \ref{prop2}, we know that the modulus of $A$ exists and $|A|$ is a Riesz multimorphism. By \cite[Proposition 3.1]{geraldoluis}, the positive $m$-linear operator $AR_{m}^{\rho}(|A|)$ is a Riesz multimorphism on $(E_{1}^{\sim})_{n}^{\sim}\times\cdots\times (E_{m}^{\sim})_{n}^{\sim}$, therefore, $AR_{m}^{\rho}(|A|)_{x_{1}^{\sim\sim},\stackrel{\hat{i}}{\ldots},x_{m}^{\sim\sim}} $ 
    is a Riesz homomorphism. By \cite[Theorem 2.14]{alip},
   \begin{equation} AR_{m}^{\rho}(|A|)_{x_{1}^{\sim\sim},\stackrel{\hat{i}}{\ldots},x_{m}^{\sim\sim}} (|x_{i}^{\sim\sim}|)\wedge AR_{m}^{\rho}(|A|)_{x_{1}^{\sim\sim},\stackrel{\hat{i}}{\ldots},x_{m}^{\sim\sim}} (|z_{i}^{\sim\sim}|)=0.\label{ybv5}
    \end{equation}
    Since    the map
    $$A \in \mathcal{L}_{r}( E_1, \ldots, E_m; F) \mapsto  AR_{m}^{\rho}(A) \in \mathcal{L}_{r}( E_1^{\sim\sim}, \ldots, E_m^{\sim\sim}; F^{\sim\sim})$$
   is a positive linear operator \cite[Theorem 2.2(c)]{geraldoluis},  $|AR_{m}^{\rho}(A)|\leq AR_{m}^{\rho}(|A|)$ because $\pm A\leq |A|$. Hence,
    \begin{align*}
      0\leq  |AR_{m}^{\rho}(A)&(x_{1}^{\sim\sim},\ldots,x_{i}^{\sim\sim},\ldots,x_{m}^{\sim\sim})|\wedge |AR_{m}^{\rho}(A)(x_{1}^{\sim\sim},\ldots,z_{i}^{\sim\sim},\ldots,x_{m}^{\sim\sim})|\\
        &\leq |AR_{m}^{\rho}(A)|(x_{1}^{\sim\sim},\ldots,|x_{i}^{\sim\sim}|,\ldots,x_{m}^{\sim\sim})\wedge |AR_{m}^{\rho}(A)|(x_{1}^{\sim\sim},\ldots,|z_{i}^{\sim\sim}|,\ldots,x_{m}^{\sim\sim})\\
        &\leq AR_{m}^{\rho}(|A|)(x_{1}^{\sim\sim},\ldots,|x_{i}^{\sim\sim}|,\ldots,x_{m}^{\sim\sim})\wedge AR_{m}^{\rho}(|A|)(x_{1}^{\sim\sim},\ldots,|z_{i}^{\sim\sim}|,\ldots,x_{m}^{\sim\sim})\\
        &=AR_{m}^{\rho}(|A|)_{x_{1}^{\sim\sim},\stackrel{\hat{i}}{\ldots},x_{m}^{\sim\sim}}(|x_{i}^{\sim\sim}|)\wedge AR_{m}^{\rho}(|A|)_{x_{1}^{\sim\sim},\stackrel{\hat{i}}{\ldots},x_{m}^{\sim\sim}}(|z_{i}^{\sim\sim}|) \stackrel{\small(\ref{ybv5})}{=} 0,
    \end{align*}
    which completes the proof.
\end{proof}

Sometimes we can reach the product of the whole biduals:

\begin{corollary}\label{p09j}
    Let $E_{1},\ldots,E_{m}, F$ be Banach lattices and let $A\colon E_1\times\cdots\times E_m\longrightarrow F$ be a regular $m$-linear operator. If the norms of $E_{1}^{\ast },\ldots, E_m^{*}$ are order continuous and $A$ is disjointness preserving, then all Arens extensions $AR_m^\rho(A) \colon E_1^{**} \times \cdots \times E_m^{**} \to F^{**}$, $\rho \in S_m$, are disjointness preserving.
\end{corollary}

\begin{proof}  By \cite[Theorem 2.4.2]{meyer} we have  $E_j^{\ast\ast}=(E_j^{\ast})_{n}^{\ast}$ for $j = 1, \ldots, m$. The result follows immediately from Theorem \ref{proporest}.
\end{proof}

Banach lattices whose duals have order continuous norms are well-known in the literature ($E^*$ has order continuous norm if and only if $E^*$ does not contain a copy of $c_0$), so we refrain from listing examples.

As to the preservation of disjointness on the whole of $E_1^{\sim\sim} \times \cdots \times E_m^{\sim\sim}$ in the setting of arbitrary Riesz spaces, we prove next a positive result for a class of regular multilinear operators, which encompasses, for instance, all regular multilinear forms.

A map taking values in a Riesz space is said to have \textit{finite lattice rank} if the sublattice generated by its range is finite dimensional. Multilinear operators of finite lattice rank proved to be useful in \cite{geraldoluis}, where the unfortunate terminology {\it finite rank} was used.

\begin{theorem} \label{teo2-6}
    Let $E_{1},\ldots,E_{m}, F$ be Riesz spaces with $F$ Archimedean and let $A\colon E_{1}\times\cdots\times E_{m}\longrightarrow F$ be a finite lattice rank regular $m$-linear  operator. If $A$  is disjointness preserving, then all Arens extensions $AR_{m}^{\rho}(A) \colon E_{1}^{\sim\sim} \times \cdots \times E_{m}^{\sim\sim} \to F^{\sim\sim}$ of $A$, $\rho\in S_{m}$, are disjointness preserving.
\end{theorem}

\begin{proof}
   By Proposition \ref{prop2} the modulus of $A$ exists, is a Riesz multimorphism and satisfies
    $$|A|(|x_1|, \ldots, |x_m|) = |A(x_1, \ldots, x_m)| \mbox{ for all } x_1 \in E_1, \ldots, x_m \in E_m. $$
    Since $A$ has finite lattice rank, its range is contained in a finite dimensional sublattice $G$ of $F$. Applying the equality above and using that $G$ is a vector sublattice of $F$,
  $$|A|(x_1, \ldots, x_m) = |A(x_1, \ldots, x_m)| \in G \mbox{ for all } x_1 \in E_1^+, \ldots, x_m \in E_m^+. $$
 Given $x_1 \in E_1, \ldots, x_m \in E_m$, writing $x_j = x_j^+ -x_j^-$ for each $j$, the $m$-linearity of $|A|$ and the sentence above give $|A|(x_1, \ldots,x_m) \in G$, proving that $|A|$ has finite lattice rank. So, $|A|$ is a Riesz multimorphism of finite lattice rank, therefore every Arens extension $AR_{m}^{\rho}(|A|), \rho\in S_{m}$, of $|A|$ is also a Riesz multimorphism by \cite[Theorem 3.2]{geraldoluis}. Now the proof goes along, {\it verbatim}, the same steps of the proof of Theorem \ref{proporest}.
\end{proof}

Although this can be checked directly, it is interesting to mention that Theorem \ref{teo2-6} implies that the bilinear operator $A$ from the counterexample given in Section 4 fails to have finite lattice rank. 

We will now move in the direction of finding situations
where each Arens extension of a disjointness preserving $m$-linear operator preserves disjointness in $j$-th coordinate on the whole of $E_j^{\sim\sim}$ for some $j \in \{1,\ldots,m\}$, whenever the remaining  coordinates are fixed at $(E_{i}^{\sim})_{n}^{\sim}, i=1,\ldots,m, i\neq j$. Two lemmas are needed. The symbol $\stackrel{\hat{j}}{\ldots}$ means that the $j$-th coordinate has been omitted.

\begin{lemma}\label{lemma5}
    Let $E_{1},\ldots,E_{m}, F$ be Riesz spaces with $F$ 
    Archimedean, and let $A\colon E_{1}\times\cdots\times E_{m}\longrightarrow F$  be a regular $m$-linear operator. If $A$ is disjointness preserving, $\rho\in S_{m}$, $j\in\{1,\ldots,m\}$, and $x_{i}\in E_{i}, i=1,\ldots,m$, $i\neq j$, then the linear operator $$AR_{m}^{\rho}(A)_{J_{E_1}(x_1),\stackrel{\hat{j}}{\ldots},J_{E_m}(x_m)}\colon E_{j}^{\sim\sim}\longrightarrow F^{\sim\sim}$$
    is  disjointness preserving.
\end{lemma}

\begin{proof}
 Let $k\in\{1,\ldots,m\}$ be such that $\rho(k)=j$. 
    Setting $W_{l}=J_{E_{\rho(l)}}(x_{\rho(l)}), l=1,\ldots,m, l\neq k$, by the proof of \cite[Lemma 2.4]{geraldoluis} we have that, for all $x_{j}^{\sim\sim}\in E_{j}^{\sim\sim}$ and $y^{\sim}\in F^{\sim}$,
\begin{align*}
    AR_{m}^{\rho}(A)_{J_{E_1}(x_1),\stackrel{\hat{j}}{\ldots},J_{E_m}(x_m)}&(x_{j}^{\sim\sim})(y^{\sim})
    =AR_{m}^{\rho}(A)(J_{E_{1}}(x_{1}),\ldots,x_{j}^{\sim\sim},\ldots,J_{E_{m}}(x_{m}))(y^{\sim})\\
&=\big(\overline{W_{m}}\circ\cdots\circ \overline{W_{k+1}}\circ\overline{x_{j}^{\sim\sim}}\circ \overline{W_{k-1}}\circ \cdots\circ \overline{W_1}\big)((y^{\sim}\circ A)_{\rho})\\
&=x_{j}^{\sim\sim}(y^{\sim}\circ A_{x_{1},\stackrel{\hat{j}}{\ldots},x_{m}})=x_{j}^{\sim\sim}\left( \left[A_{x_{1},\stackrel{\hat{j}}{\ldots},x_{m}}\right]^{\sim}(y^{\sim})\right)\\
&= \left[A_{x_{1},\stackrel{\hat{j}}{\ldots},x_{m}}\right]^{\sim\sim}(x^{\sim\sim})(y^\sim).
\end{align*}
This proves that  $AR_{m}^{\rho}(A)_{J_{E_1}(x_1),\stackrel{\hat{j}}{\ldots},J_{E_m}(x_m)} = \left[A_{x_{1},\stackrel{\hat{j}}{\ldots},x_{m}}\right]^{\sim\sim}$.
By assumption, the regular linear operator $A_{x_{1},\stackrel{\hat{j}}{\ldots},x_{m}}$ 
is disjointness preserving, so is its second adjoint by Theorem \ref{teo1}, 
that is, $AR_{m}^{\rho}(A)_{J_{E_1}(x_1),\stackrel{\hat{J}}{\ldots},J_{E_m}(x_m)}$ is disjointness preserving.
\end{proof}

For a Riesz space $E$, recall that the weak$^{\ast}$ absolute topology $|\sigma|^{\ast}$ on $E^{\sim\sim}$ is the locally convex-solid topology generated by the family of lattice seminorms $\{p_{x^{\sim}}\colon E^{\sim \sim}\longrightarrow\mathbb{R}:x^{\sim}\in E^{\sim}\}$, where $p_{x^{\sim}}(x^{\sim\sim})=|x^{\sim\sim}|(|x^{\sim}|)$ for all $x^{\sim}\in E^{\sim}$ and $x^{\sim\sim}\in E^{\sim\sim}$. A net $(x_{\alpha}^{\sim\sim})_{\alpha}$ in $E^{\sim\sim}$ converges to $x^{\sim\sim}\in E^{\sim\sim}$ in the weak$^{\ast}$ absolute topology  if and only if $|x_{\alpha}^{\sim\sim}-x^{\sim\sim}|(|x^{\sim}|)\longrightarrow 0$ for every $x^{\sim}\in E^{\sim}$. By \cite[Theorem 1.7(6) and p.\,8]{alip}, the lattice operations are $|\sigma|^{\ast}$-$|\sigma|^{\ast}$ continuous.

Next lemma can be regarded as an echo, in the lattice context, of the Davie-Gamelin formula for Arens extensions in Banach spaces (see, e.g., \cite[p.\,412]{din}).

\begin{lemma}\label{lemma6}
   Let $A\colon E_{1}\times\cdots\times E_{m}\longrightarrow F$  be a regular $m$-linear operator between Riesz spaces, let $\rho\in S_{m}$ be given, and let $j \in \{1,\ldots, m\}$ be such that $j=\rho(1)$. Given $x_i^{\sim\sim} \in E_i ^{\sim\sim}, i = 1,\ldots, m$, suppose that, for each $i \neq j$, there exists a net $(x_{\alpha_{i}})_{\alpha_{i}}$ in $E_{i}$ such that $x_{i}^{\sim\sim}=|\sigma|^{\ast}$-$\displaystyle\lim_{\alpha_{i}}J_{E_{i}}(x_{\alpha_{i}})$ in $E_{i}^{\sim\sim}$.  Then, 
    $$AR_{m}^{\rho}(A)(x_{1}^{\sim\sim},\ldots,x_{j}^{\sim\sim},\ldots,x_{m}^{\sim\sim})=
    |\sigma|^{\ast}\mbox{-}\lim_{\alpha_{\rho(m)}}\cdots|\sigma|^{\ast}\mbox{-}\lim_{\alpha_{\rho(2)}} AR_{m}^{\rho}(A)_{J_{E_1}(x_{\alpha_1}),\stackrel{\hat{j}}{\cdots},J_{E_m}(x_{\alpha_m})}(x_{j}^{\sim\sim}).$$
\end{lemma}

\begin{proof}
We prove the result for the permutation $\theta\in S_m$, see (\ref{i8m3}); the proofs for the other Arens extensions are analogous. For simplicity, we shall write $A^{*(m+1)}$ instead of $AR_m^\theta(A)$. For each $i\in \{1,\ldots,m\}$, let $x_{\theta(1)}^{\sim\sim}\in E_{\theta(1)}^{\sim\sim},\ldots,x_{\theta(i)}^{\sim\sim}\in E_{\theta(i)}^{\sim\sim},x_{\theta(i+1)}\in E_{\theta(i+1)},\ldots,x_{\theta(m)}\in E_{\theta(m)}$ be given. Since $\theta(k)=m-k+1$ for each  $k=1,\ldots,m$, by the proof of \cite[Lemma 2.4]{geraldoluis} we have
\begin{align*}
    A^{\ast (m+1)}(J_{E_{1}}(x_1),\ldots,& J_{E_{m-i}}(x_{m-i}),x_{m-i+1}^{\sim\sim},\ldots,x_{m}^{\sim\sim})(y^{\sim})\\
&=x_{\theta(i)}^{\sim\sim}\Big(\Big((\overline{x_{\theta(i-1)}^{\sim\sim}}\circ\cdots\circ \overline{x_{\theta(1)}^{\sim\sim}})((y^{\sim}\circ A)_{\theta})\Big)^{i}(x_{m-i},\ldots,x_{1})\Big)
\end{align*}
for every $y^{\sim}\in F^{\sim}$. Let us prove that the map
$$x_{\theta(i)}^{\sim\sim} \in E_{\theta(i)}^{\sim\sim} \mapsto A^{\ast (m+1)}(J_{E_{1}}(x_1),\ldots,J_{E_{m-i}}(x_{m-i}), x_{m-i+1}^{\sim\sim},\ldots,x_{m}^{\sim\sim}) \in F^{\sim\sim}  $$
is $|\sigma|^*$-$|\sigma|^*$ continuous.
To do so, let $(x_{\alpha_{\theta(i)}}^{\sim\sim})_{\alpha_{\theta(i)}}$ be a net in $E_{\theta(i)}^{\sim\sim}$ such that $x_{\alpha_{\theta(i)}}^{\sim\sim}\xrightarrow{\,\, |\sigma|^{\ast} } x_{\theta(i)}^{\sim\sim}$. Since $A$ is regular, there are  positive $m$-linear operators $A_{1}, A_{2}\colon E_{1}\times\cdots\times E_{m}\longrightarrow F$ such that $A=A_1-A_2$, 
thus
$$ |A^{*(m+1)}| = |A_1^{*(m+1)} - A_2^{*(m+1)}| \leq A_1^{*(m+1)} + A_2^{*(m+1)} =: B^{*(m+1)},$$
For every $y^{\sim}\in F^{\sim}$,
\begin{align*}
    |  A^{\ast (m+1)}&(J_{E_{1}}(x_1),\ldots,J_{E_{m-i}}(x_{m-i}), x_{\alpha_{m-i+1}}^{\sim\sim}-x_{m-i+1}^{\sim\sim},\ldots,x_{m}^{\sim\sim})|(|y^{\sim}|)\\
    &\leq |A^{\ast (m+1)}|(J_{E_{1}}(|x_1|),\ldots,J_{E_{m-i}}(|x_{m-i}|), |x_{\alpha_{m-i+1}}^{\sim\sim}-x_{m-i+1}^{\sim\sim}|,\ldots,|x_{m}^{\sim\sim}|)(|y^{\sim}|)\\
    &\leq B^{\ast (m+1)}(J_{E_{1}}(|x_1|),\ldots,J_{E_{m-i}}(|x_{m-i}|), |x_{\alpha_{m-i+1}}^{\sim\sim}-x_{m-i+1}^{\sim\sim}|,\ldots,|x_{m}^{\sim\sim}|)(|y^{\sim}|)\\
&=|x_{\alpha_{\theta(i)}}^{\sim\sim}-x_{\theta(i)}^{\sim\sim}|\Big(\Big((\overline{|x_{\theta(i-1)}^{\sim\sim}|}\circ\cdots\circ \overline{|x_{\theta(1)}^{\sim\sim}|})((|y^{\sim}|\circ B)_{\theta})\Big)^{i}(|x_{m-i}|,\ldots,|x_{1}|)\Big)\longrightarrow 0,
\end{align*}
from which it follows that the net $(A^{\ast (m+1)}(J_{E_{1}}(x_1),\ldots,J_{E_{m-i}}(x_{m-i}), x_{\alpha_{m-i+1}}^{\sim\sim},\ldots,x_{m}^{\sim\sim}))_{\alpha_{\theta(i)}}$  converges to $A^{\ast (m+1)}(J_{E_{1}}(x_1),\ldots,J_{E_{m-i}}(x_{m-i}), x_{m-i+1}^{\sim\sim},\ldots,x_{m}^{\sim\sim})$ in the weak$^{\ast}$ absolute topology $|\sigma|^{\ast}$ of $F^{\sim\sim}$. The desired $|\sigma|^*$-$|\sigma|^*$ continuity has been proved.

Recall that, for the permutation we are working with, it holds $j=\theta(1)=m$. By assumption, for each $k = 1, \ldots, m-1$, there exists a net 
$(x_{\alpha_{k}})_{\alpha_{k}}$ in $E_{k}$ 
 such that $J_{E_{k}}(x_{\alpha_{k}})\xrightarrow{\,\, |\sigma|^{\ast} } x_{k}^{\sim\sim}$. 
 Applying the $|\sigma|^*$-$|\sigma|^*$ continuity proved above with $i=2$, we get
 \begin{equation}\label{parai2}
     A^{\ast (m+1)}(J_{E_{1}}(x_1),\ldots,J_{E_{m-2}}(x_{m-2}),x_{m-1}^{\sim\sim},x_{m}^{\sim\sim})=|\sigma|^{\ast}\mbox{-}\displaystyle\lim_{\alpha_{m-1}} A^{\ast (m+1)}_{x_1,\ldots,x_{\alpha_{m-1}}}(x_{m}^{\sim\sim}),
 \end{equation}
 for all $x_{1}\in E_1,\ldots,x_{m-2}\in E_{m-2}$. 
  Repeating the procedure with $i=3$ and using (\ref{parai2}), we have
 \begin{align*}
     A^{\ast (m+1)}(J_{E_{1}}(x_1),\ldots,&J_{E_{m-3}}(x_{m-3}),x_{m-2}^{\sim\sim},x_{m-1}^{\sim\sim},x_{m}^{\sim\sim})\\
     &=|\sigma|^{\ast}\mbox{-}\displaystyle\lim_{\alpha_{m-2}}A^{\ast (m+1)}(J_{E_{1}}(x_1),\ldots,J_{E_{m-2}}(x_{\alpha_{m-2}}),x_{m-1}^{\sim\sim},x_{m}^{\sim\sim})\\
     &=|\sigma|^{\ast}\mbox{-}\displaystyle\lim_{\alpha_{m-2}}|\sigma|^{\ast}-\displaystyle\lim_{\alpha_{m-1}} A^{\ast (m+1)}_{x_1,\ldots,x_{\alpha_{m-2}},x_{\alpha_{m-1}}}(x_{m}^{\sim\sim}),
 \end{align*}
for all $x_{1}\in E_{1},\ldots,x_{m-3}\in E_{m-3}, x_{m-1}^{\sim\sim}\in E_{m-1}^{\sim\sim}$. 
Proceeding as above $(m-3)$-times, it follows that
$$A^{\ast (m+1)}(x_{1}^{\sim\sim},\ldots,x_{m-1}^{\sim\sim},x_{m}^{\sim\sim})=|\sigma|^{\ast}\mbox{-}\displaystyle\lim_{\alpha_{1}}|\cdots\sigma|^{\ast}\mbox{-}\displaystyle\lim_{\alpha_{m-1}} A^{\ast (m+1)}_{x_{\alpha_{1}},\ldots,x_{\alpha_{m-2}},x_{\alpha_{m-1}}}(x_{m}^{\sim\sim}).$$
\end{proof}

Next we shall use that, for two vectors $x,y$ in a Riesz space,
\begin{equation} \label{rt9u}x \perp y \Longleftrightarrow |x+y| = |x-y|. \end{equation}

\begin{proposition}\label{propo3}
    Let $A\colon E_1\times\cdots\times E_m\longrightarrow F$ be a regular $m$-linear operator between Archimedean Riesz spaces and let $\rho\in S_m$ be given. If $A$ is disjointness preserving, then the Arens extension $AR_{m}^{\rho}(A)$ of $A$ is disjointness preserving in the $\rho(1)$-th coordinate on $(E_{1}^{\sim})_{n}^{\sim}\times\cdots\times(E_{\rho(1)-1}^{\sim})_{n}^{\sim}\times E_{\rho(1)}^{\sim\sim}\times (E_{\rho(1)+1}^{\sim})_{n}^{\sim}\times\cdots\times (E_{m}^{\sim})_{n}^{\sim}$.
\end{proposition}

\begin{proof} For simplicity we write 
$j=\rho(1)$. Let $x_{i}^{\sim\sim}\in (E_{i}^{\sim})_{n}^{\sim}, i=1,\ldots,m, i\neq j$, be given. For every such $i$, by \cite[Theorem 2]{grobler} there exists a net $(x_{\alpha_{i}})_{\alpha_{i}}$ in $E_{i}$ such that $J_{E_{i}}(x_{\alpha_{i}})\xrightarrow{\,\, |\sigma|^{\ast} } x_{i}^{\sim\sim}$. Let $x_{j}^{\sim\sim}, z_{j}^{\sim\sim}\in E_{j}^{\sim\sim}$ be such that $x_{j}^{\sim\sim}\perp z_{j}^{\sim\sim}$. 
For all $x_{1}\in E_1, \ldots,x_{j-1}\in E_{j-1},x_{j+1}\in E_{j+1},\ldots,x_{m}\in E_{m}$, Lemma \ref{lemma5} gives
$$AR_{m}^{\rho}(A)_{J_{E_1}(x_1),\stackrel{\hat{j}}{\ldots},J_{E_m}(x_m)}(x_j^{\sim\sim}) \perp AR_{m}^{\rho}(A)_{J_{E_1}(x_1),\stackrel{\hat{j}}{\ldots},J_{E_m}(x_m)}(z_j^{\sim\sim}).  $$
Combining (\ref{rt9u}) with the linearity of $AR_{m}^{\rho}(A)$ in the $j$-th coordinate, it follows that
    \begin{equation}\label{eq}
        |AR_{m}^{\rho}(A)_{J_{E_1}(x_1),\stackrel{\hat{j}}{\ldots},J_{E_m}(x_m)}(x_{j}^{\sim\sim}+z_{j}^{\sim\sim})|
        =|AR_{m}^{\rho}(A)_{J_{E_1}(x_1),\stackrel{\hat{j}}{\ldots},J_{E_m}(x_m)}(x_{j}^{\sim\sim}-z_{j}^{\sim\sim})|.
    \end{equation}
Using the $|\sigma|^{\ast}\mbox{-}|\sigma|^{\ast}$ continuity of the lattice operations together with Lemma \ref{lemma6} and (\ref{eq}), we have
    \begin{align*}
        |AR_{m}^{\rho}(A)(x_{1}^{\sim\sim},\ldots,&x_{j}^{\sim\sim}+z_{j}^{\sim\sim},\ldots,x_{m}^{\sim\sim})|\\
&=|\sigma|^{\ast}\mbox{-}\lim_{\alpha_{\rho(m)}}\cdots|\sigma|^{\ast}\mbox{-}\lim_{\alpha_{\rho(2)}} |AR_{m}^{\rho}(A)_{J_{E_1}(x_{\alpha_1}),\stackrel{\hat{j}}{\ldots},J_{E_m}(x_{\alpha_m})}(x_{j}^{\sim\sim}+z_{j}^{\sim\sim})|\\
&=|\sigma|^{\ast}\mbox{-}\lim_{\alpha_{\rho(m)}}\cdots|\sigma|^{\ast}\mbox{-}\lim_{\alpha_{\rho(2)}} |AR_{m}^{\rho}(A)_{J_{E_1}(x_{\alpha_1}),\stackrel{\hat{j}}{\ldots},J_{E_m}(x_{\alpha_m})}(x_{j}^{\sim\sim}-z_{j}^{\sim\sim})|\\
 &= |AR_{m}^{\rho}(A)(x_{1}^{\sim\sim},\ldots,x_{j}^{\sim\sim}-z_{j}^{\sim\sim},\ldots,x_{m}^{\sim\sim})|.
    \end{align*}
   The result follows by calling on (\ref{rt9u}) once again. 
\end{proof}

Our next purpose is to show that, for operators between Banach lattices with suitable properties, we can get disjointness preservation on the product of the whole of the biduals. For simplicity, we shall deal with the bilinear case. We will use that every Banach lattice is Archimedean \cite[Proposition 1.1.8(i)]{meyer}.

 Recall that $\theta$ is the permutation of $\{1,2\}$ given by $\theta(1) = 2, \theta(2)= 1$. The other permutation, that is, the identity, shall be denoted by ${\rm id}$.

 \begin{theorem}\label{iu8h} Let $E_1, E_2,F$ be Banach lattices such that the lattice operations in $F^{**}$ are weak*-sequentially continuous. Let $A\colon E_{1}\times E_{2}\longrightarrow F$  be a disjointness preserving regular bilinear operator.\\
 {\rm (a)} If $E_1^*$ is separable, then, for every $ x_{1}^{\ast\ast}\in E_{1}^{\ast\ast}$, the operator $AR_{2}^{\theta}(A)_{x_1^{**}} \colon E_2^{**} \to F^{**}$ is disjointness preserving. \\
 {\rm (b)} If $E_2^*$ is separable, then, for every $ x_{2}^{\ast\ast}\in E_{2}^{\ast\ast}$, the operator $AR_{2}^{\rm id}(A)_{x_2^{**}} \colon E_1^{**} \to F^{**}$ is disjointness preserving.
 \end{theorem}

 \begin{proof} (a) Let $x_{2}^{\ast\ast}, z_{2}^{\ast\ast}\in E_{2}^{\ast\ast}$ be such that $x_{2}^{\ast\ast} \perp z_{2}^{\ast\ast}$. By 
 Lemma \ref{lemma5} we have $$AR_{2}^{\theta}(A)(J_{E_{1}}(x_{1}), x_{2}^{\ast\ast})\perp AR_{2}^{\theta}(A)(J_{E_{1}}(x_{1}), z_{2}^{\ast\ast}) \text{ for every } x_{1}\in E_{1}.$$
Let $x_1^{**} \in E_1^{**} $ be such that $\norma{x_1^{**}} = 1$. Since $E_1^{*}$ is separable, the weak$^{\ast}$ topology of $E_1^{**}$ is metrizable on bounded sets of $E_1^{**}$. So, as $J_{E_1}(B_{E_{1}})$ is weak$^{\ast}$-dense in $B_{E_1^{**}}$, there exists a sequence $(x_{1, n})_n$ in $B_{E_1}$ such that $J_{E_1}(x_{1,n}) \xrightarrow{\,\omega^{\ast}}x_1^{**}$ in $E_{1}^{**}$.
Since the Arens extension $AR_{2}^{\theta}(A)$ is weak$^{\ast}$-weak$^*$ continuous in the first variable \cite[p.\,413]{din}, $AR_{2}^{\theta}(A)(x_{1}^{\ast\ast}, z_{2}^{\ast\ast})=\omega^{\ast}\mbox{\rm -}\displaystyle\lim_{n}AR_{2}^{\theta}(A)(J_{E_{1}}(x_{1,n}), z_{2}^{\ast\ast})$ and $AR_{2}^{\theta}(A)(x_{1}^{\ast\ast}, x_{2}^{\ast\ast})=\omega^{\ast}\mbox{\rm -}\displaystyle\lim_{n}AR_{2}^{\theta}(A)(J_{E_{1}}(x_{1,n}), x_{2}^{\ast\ast})$. Hen-ce,
\begin{align*}
|AR_{2}^{\theta}(A)(x_{1}^{\ast\ast},& x_{2}^{\ast\ast})|\wedge |AR_{2}^{\theta}(A)(x_{1}^{\ast\ast}, z_{2}^{\ast\ast})|\\
&=\omega^{\ast}\mbox{-}\lim_{n}|AR_{2}^{\theta}(A)(J_{E_{1}}(x_{1,n}), x_{2}^{\ast\ast})|\wedge \omega^{\ast}\mbox{-}\lim_{n}|AR_{2}^{\theta}(A)(J_{E_{1}}(x_{1,n}), z_{2}^{\ast\ast})|\\
&=\omega^{\ast}\mbox{-}\lim_{n}\left(|AR_{2}^{\theta}(A)(J_{E_{1}}(x_{1,n}), x_{2}^{\ast\ast})|\wedge|AR_{2}^{\theta}(A)(J_{E_{1}}(x_{1,n}), z_{2}^{\ast\ast})|\right)=0.
\end{align*}
Thus far we have $AR_{2}^{\theta}(A)(x_{1}^{\ast\ast}, x_{2}^{\ast\ast}) \perp AR_{2}^{\theta}(A)(x_{1}^{\ast\ast}, z_{2}^{\ast\ast})$ for every $x_1^{**} \in E_1^{**} $ with $\norma{x_1^{**}} = 1$. For $x_1^{**} = 0$ this is obvious, and for arbitrary nonzero $x_1^{**}$ it is enough to normalize $x_1^{**}$ and apply the bilinearity of $AR_{2}^{\theta}(A)$ together with \cite[Ex.\,2(a),\,p.\,21]{alip}.  

(b) The proof is analogous.
 \end{proof}

 Recall that a Banach space $E$ is {\it Arens regular} if every bounded linear operator
from $E$ to $E^{\ast}$ is weakly compact. Examples shall be provided soon. 

\begin{corollary} \label{cor4}
    Let $E$ be an Arens regular Banach lattice with separable dual, and let $F$ be a Banach lattice such that the lattice operations in $F^{**}$ are weak*-sequentially continuous. Then, the unique Arens extension of any disjointness preserving regular bilinear operator from $E^2$ to $F$ is disjointness preserving on $(E^{**})^2$.
\end{corollary}

\begin{proof} By \cite[Proposition 6.14]{din}, the two Arens extensions of every continuous bilinear operator from $E^2$ to any Banach space are coincident. The result follows from Theorem \ref{iu8h}.
\end{proof}

Next, we outline cases where Corollary \ref{cor4} can be directly applied. 

\begin{example}  \rm
{\rm (a)} Of course, only nonreflexive domain spaces are of interest. 
Examples of nonreflexive Arens regular Banach lattices with separable duals are the following: (i) $c_0$ and $c$. (ii) $c_0(\{\ell_1\}_k)$, which is the $c_0$-sum of the sequence $(\ell_1^k)_{k \in \mathbb{N}}$, where $\ell_1^k = (\mathbb{R}^k, \|\cdot\|_1)$. Of course, this is a nonreflexive Banach lattice with separable dual. For its Arens regularity, see \cite[Remark 1.4(d)]{transac}. (iii) Infinite dimensional separable $AM$-spaces with unit and with separable dual, that is,  $C(K)$-spaces, where $K$ is a countable compact Hausdorff space. Its dual $C(K)^*$ is an $AL$-space, which is an $L_1(\mu)$-space. By a result due to Grothendieck (see \cite[p.\,251]{pel}), every operator from $C(K)$ to $C(K)^*$ factors through a Hilbert space, hence is weakly compact. 
\\(b) By \cite[Proposition 2.1]{boydmiranda}, any Banach lattice whose dual has order continuous norm is Arens regular, so any nonreflexive separable Banach lattice whose norm is order continuous fits  our purpose. 
Examples: $c_0(E)$ for any separable reflexive Banach lattice $E$ and 
$\ell_q(c_0),1 < q < \infty$.
\\{\rm (c)}  For each Banach lattice $F$ such that $F^*$ is discrete and order continuous, the lattice operations in $F^{**}$ are weak*-sequentially continuous \cite[Theorem 6.6]{wnukbook}. For instance, $c_0(\Gamma)$ and $\ell_p(\Gamma)$, $1 < p < \infty$, are examples of such Banach lattices for every index set $\Gamma$.
\end{example}

\begin{remark}\rm
   Theorem \ref{iu8h} and Corollary \ref{cor4} can be proved, {\it mutatis mutandis}, for disjointness preserving regular $m$-linear operators, $m \geq 2$. The proofs are analogous using that every Arens extension of a continuous $m$-linear operator is weak$^{\ast}$-weak$^{\ast}$ continuous in one variable (see the proof of \cite[Lemma 2.4]{geraldoluis}). For the uniqueness of the extension in Corollary \ref{cor4}, see \cite[Theorem 1]{bombal}.
\end{remark}

 We finish the paper by finding conditions on the Banach lattice $F$ so that all Arens extensions of every disjointness preserving regular $F$-valued $m$-linear operator are disjointness preserving on the product of the whole biduals. 

\begin{lemma} \label{prop3}
    Let $E_1, \dots, E_m, \, F$ be Riesz spaces
    , let
    $A \colon E_1 \times \cdots \times E_m \to F$ be a regular disjointness preserving $m$-linear operator, and let $\rho \in S_m$ be given.
    If $y^{\sim} \in F^\sim$ is disjointness preserving, then
    $J_{F^\sim}(y^{\sim}) \circ AR_{m}^{\rho}(A)$ is disjointness preserving on $E_1^{\sim\sim} \times \cdots \times E_m^{\sim\sim}$ and
     \begin{align*}|AR_{m}^{\rho}(A)(x_1^{\sim\sim}, \dots, x_m^{\sim\sim})|(|y^{\sim}|) &= |AR_{m}^{\rho}(A)(|x_1^{\sim\sim}|, \dots, |x_m^{\sim\sim}|)(y^{\sim})| \\
    &=|AR_{m}^{\rho}(A)(x_1^{\sim\sim}, \dots, x_m^{\sim\sim})(y^{\sim})|
    \end{align*}
    for all $x_1^{\sim\sim} \in E_1^{\sim \sim}, \dots, x_m^{\sim\sim} \in E_m^{\sim \sim}$.
\end{lemma}

\begin{proof}
     Since $A$ is regular and disjointness preserving and $y^{\sim} \in F^\sim$ is disjointness preserving, $y^{\sim} \circ A \colon E_1 \times \cdots \times E_m \to \R$ is a regular disjointness preserving $m$-linear operator. 
     Applying Theorem \ref{teo2-6}   we get that $AR_{m}^{\rho}(y^{\sim} \circ A)$ is also disjointness preserving. The proof of \cite[Proposition 3.4]{geraldoluis} yields that $AR_{m}^{\rho}(y^{\sim} \circ A) = J_{F^\sim}(y^{\sim}) \circ AR_{m}^{\rho}(A)$ holds, proving that $J_{F^\sim}(y^{\sim}) \circ AR_{m}^{\rho}(A)$ is disjointness preserving and that
    \begin{equation} \label{eq2}
        |J_{F^\sim}(y^{\sim}) \circ AR_{m}^{\rho}(A)(|x_1^{\sim\sim}|, \dots, |x_m^{\sim\sim}|)| = |J_{F^\sim}(y^{\sim}) \circ AR_{m}^{\rho}(A)(x_1^{\sim\sim}, \dots, x_m^{\sim\sim})|
    \end{equation}
    holds for all $x_1^{\sim\sim} \in E_1^{\sim \sim}, \dots, x_m^{\sim\sim} \in E_m^{\sim \sim}$. Observe also that $J_{F^\sim}(y^{\sim})$ is disjointness preserving by Theorem \ref{teo1}(1), hence
    \begin{equation} \label{eq3}
         |J_{F^\sim}(y^{\sim})|(|y^{\sim\sim}|) = |J_{F^\sim}(y^{\sim})(y^{\sim\sim})|
    \end{equation}
    for every $y^{\sim\sim} \in F^{\sim \sim}$. Combining (\ref{eq2}) and (\ref{eq3}), we have         \begin{align*}
        |AR_{m}^{\rho}(A)(|x_1^{\sim\sim}|, \dots, |x_m^{\sim\sim}|)(y^{\sim})| & = |(J_{F^\sim}(y^{\sim}) \circ AR_{m}^{\rho}(A))(|x_1^{\sim\sim}|, \dots, |x_m^{\sim\sim}|)| \\
        & = |(J_{F^\sim}(y^{\sim}) \circ AR_{m}^{\rho}(A))(x_1^{\sim\sim}, \dots, x_m^{\sim\sim})| \\
        & = |AR_{m}^{\rho}(A)(x_1^{\sim\sim}, \dots, x_m^{\sim\sim})(y^{\sim})| \\
        & \leq |AR_{m}^{\rho}(A)(x_1^{\sim\sim}, \dots, x_m^{\sim\sim})| (|y^{\sim}|)  \\
        & = J_{F^\sim}(|y^{\sim}|)( |AR_{m}^{\rho}(A)(x_1^{\sim\sim}, \dots, x_m^{\sim\sim})|) \\
        & = |J_{F^\sim}(y^{\sim})| ( |AR_{m}^{\rho}(A)(x_1^{\sim\sim}, \dots, x_m^{\sim\sim})|) \\
        & = |J_{F^\sim}(y^{\sim})(AR_{m}^{\rho}(A)(x_1^{\sim\sim}, \dots, x_m^{\sim\sim}))| \\
        & = |AR_{m}^{\rho}(A)(|x_1^{\sim\sim}|,\dots, |x_m^{\sim\sim}|) (y^{\sim}) |
    \end{align*}
    for all $x_1^{\sim\sim} \in E_1^{\sim \sim}, \dots, x_m^{\sim\sim} \in E_m^{\sim \sim}$.
\end{proof}


 Next result can be regarded as a version of \cite[Proposition 3.6]{geraldoluis} for disjointness preserving multilinear operators.

\begin{theorem} \label{prop4}
Let $E_1, \ldots, E_m, F$ be Banach lattices, 
let $A \in {\cal L}_r(E_1, \ldots, E_m;F)$ be a disjointness preserving  operator, and let $\rho \in S_m$ be given. Then:\\
{\rm (a)} $y^{***} \circ AR_m^{\rho}(A)$ is a disjointness preserving $m$-linear operator for every weak*-continuous disjointness preserving functional  $y^{***} \in F^{***}$.\\
{\rm (b)} Fix $i\in\{1,\ldots,m\}$ and let $w_{i}^{\ast\ast}, z_{i}^{\ast\ast} \in E_i^{\ast \ast}$ be such that
$w_{i}^{\ast\ast}\perp z_{i}^{\ast\ast}$. Then
$$\big(|AR_{m}^{\rho}(A)(x_{1}^{**},\ldots,w_{i}^{\ast\ast},\ldots,x_{m}^{**})|\wedge |AR_{m}^{\rho}(A)(x_{1}^{**},\ldots,z_{i}^{\ast\ast},\ldots,x_{m}^{**})|\big)(y^{*})=0$$
for all $x_{j}^{**}\in E_{j}^{**}, j\neq i$ and any $y^* \in \overline{{\rm span}}\{\varphi\in F^{\ast}:\varphi \text{ is disjointness preserving}\}.$
\end{theorem}
\begin{proof}
{\rm (a)} Let $y^{***} \in F^{***}$ be a weak*-continuous disjointness preserving functional, and let 
$y^{\ast}\in F^{\ast}$ be such that $y^{\ast\ast\ast}=J_{F^{\ast}}(y^{\ast})$. As $y^{\ast\ast\ast}$ is disjointness preserving
and $J_{F}$ is a Riesz homomorphisms, $y^{\ast}=y^{\ast\ast\ast}\circ J_{F}$ is disjointness preserving. Now, Lemma \ref{prop3} yields that $y^{***} \circ AR_m^{\rho}(A) = J_{F^{\ast}}(y^{\ast})\circ AR_{m}^{\rho}(A)$ is a disjointness preserving operator.\\
(b) For each $j \neq i$, take
$x_{j}^{\ast\ast}\in E_{j}^{\ast\ast}$.  Given $y^{\ast}\in \overline{{\rm span}}\{\varphi\in F^{\ast}:\varphi \text{ is a disjointness preserving}\}$, there exists a sequence $(y_{n}^{\ast})_{n=1}^{\infty}$ in ${\rm span}\{\varphi\in F^{\ast}:\varphi \text{ is a disjointness preserving}\}$  such that $y_{n}^{\ast}\longrightarrow y^{\ast}$ in $F^*$, hence $|y_{n}^{\ast}|\longrightarrow |y^{\ast}|$. Moreover, for each  $n\in\mathbb{N}$ there exist disjointness preserving linear functionals $\varphi^1_n,\ldots,\varphi_n^{k_n}\in F^{\ast}$ and scalars $\alpha_n^1,\ldots,\alpha_n^{k_n} $ such that  $y_{n}^{\ast}=\sum\limits_{j=1}^{k_n}\alpha_n^{j}\varphi_n^j$. For each $t_i^{**} \in  E_{i}^{\ast \ast}$, write
 $$\displaystyle AR_{m}^{\rho}(A)_{i}(t_i^{**}):=AR_{m}^{\rho}(A)(x_{1}^{\ast\ast},\ldots,x_{i-1}^{\ast\ast},t_i^{**}, x_{i+1}^{\ast\ast},\ldots,x_{m}^{\ast\ast}).$$
It follows from Lemma \ref{prop3} that
$ \displaystyle |AR_{m}^{\rho}(A)_{i}(t_i^{**})|(|\varphi|) =  |AR_{m}^{\rho}(A)_{i}(t_i^{**})(\varphi)|$ for every disjointness preserving $\varphi \in F^\ast$ and every $t_i^{**} \in E_i^{\ast\ast}$. Therefore,
\begin{align*}
  { |\big(|AR_{m}^{\rho}}&{(A)_{i} (w_{i}^{\ast\ast})|  \wedge |AR_{m}^{\rho}(A)_{i}(z_{i}^{\ast\ast})|\big)(y^{*})|\leq \big(|AR_{m}^{\rho}(A)_{i} (w_{i}^{\ast\ast})|  \wedge |AR_{m}^{\rho}(A)_{i}(z_{i}^{\ast\ast})|\big)(|y^{*}|)} \\
    &= \lim_{n\rightarrow \infty} \big(|AR_{m}^{\rho}(A)_{i}(w_{i}^{\ast\ast})|\wedge |AR_{m}^{\rho}(A)_{i}(z_{i}^{\ast\ast})|\big)(|y_{n}^{\ast}|)\\
    &=\lim_{n\rightarrow \infty} \big(|AR_{m}^{\rho}(A)_{i}(w_{i}^{\ast\ast})|\wedge |AR_{m}^{\rho}(A)_{i}(z_{i}^{\ast\ast})|\big)\Big(\Big|\sum\limits_{j=1}^{k_n}\alpha_n^{j}\varphi_n^j\Big|\Big)\\
    &\leq \lim_{n\rightarrow \infty} \displaystyle\sum_{j=1}^{k_n} |\alpha_n^{j}| \big(|AR_{m}^{\rho}(A)_{i}(w_{i}^{\ast\ast})|\wedge |AR_{m}^{\rho}(A)_{i}(z_{i}^{\ast\ast})|\big)(|\varphi_n^{j}|)\\
    &\leq  \lim_{n\rightarrow \infty} \displaystyle\sum_{j=1}^{k_n} |\alpha_n^{j}| \big(|AR_{m}^{\rho}(A)_{i}(w_{i}^{\ast\ast})|(|\varphi_n^{j}|)\wedge |AR_{m}^{\rho}(A)_{i}(z_{i}^{\ast\ast})|(|\varphi_n^{j}|)\big)\\
    & =  \lim_{n\rightarrow \infty} \displaystyle\sum_{j=1}^{k_n} |\alpha_n^{j}| \big(|AR_{m}^{\rho}(A)_{i}(w_{i}^{\ast\ast})(\varphi_n^{j})|\wedge |AR_{m}^{\rho}(A)_{i}(z_{i}^{\ast\ast})(\varphi_n^{j})|\big)\\
    &=   \lim_{n\rightarrow \infty} \displaystyle\sum_{j=1}^{k_n} |\alpha_n^{j}| \big(|(J_{F^{\ast}}(\varphi_n^{j})\circ AR_{m}^{\rho}(A)_{i})(w_{i}^{\ast\ast})|\wedge |(J_{F^{\ast}}(\varphi_n^{j})\circ AR_{m}^{\rho}(A)_{i})(z_{i}^{\ast\ast})|\big)=0,
\end{align*}
where the last equality follows from Lemma \ref{prop3}.
\end{proof}

Our last main result is a straightforward consequence of Theorem \ref{prop4}(b).

\begin{corollary} \label{cor2}
    Let $F$ be a Banach lattice such that $F^\ast$ has a Schauder basis consisting of disjointness preserving linear functionals. Then
    all Arens extensions of any $F$-valued disjointness preserving regular multilinear operator are disjointness preserving on the product of the whole biduals.
\end{corollary}

Next, we list some cases where Corollary \ref{cor2} applies.

\begin{example}\rm Since Riesz homomorphisms are disjointness preserving, from \cite[Example 3.9]{geraldoluis} it follows that Corollary \ref{cor2} applies to the following Banach lattices $F$: $c_0, \ell_p, 1<p < \infty$, any Banach lattice whose order is given by an 1-unconditional Schauder basis, Tsirelson's original space $T^*$ and its dual $T$, Schreier's space $S$, and the predual $d_*(w,1)$ of the Lorentz sequence space $d(w,1)$.
\end{example}

 We remark that the assumptions of Theorem \ref{teo2-6} and Corollary \ref{cor2} are independent. On the one hand, Theorem \ref{teo2-6} holds for Riesz spaces, whereas Corollary \ref{cor2} holds only for Banach lattices. In particular, the assumptions of Theorem \ref{teo2-6} do not imply those of Corollary \ref{cor2}. On the other hand, the dual of $c_{0}$ has a Schauder basis consisting of disjointness preserving linear functionals. Now, consider the bilinear operator $A \colon c_{0} \times c_{0} \longrightarrow c_{0}$ defined by $A(x,y) = (x_{n}y_{n})_{n=1}^{\infty}$. Clearly, $A$ is disjointness preserving, since it is a $2$-morphism, but it does not have finite lattice rank, because the canonical sequence $(e_{n})_{n=1}^{\infty}$ is contained in $A(c_{0} \times c_{0})$. Hence, the assumptions of Corollary \ref{cor2} do not imply those of Theorem \ref{teo2-6}.












\paragraph*{Declaration on the use of generative AI.} The original research questions, conceptual framework, and main results developed in this paper are due to the authors. OpenAI's ChatGPT was used exclusively to obtain the counterexample in Section \ref{contra}. The resulting mathematical arguments were subsequently carefully checked, corrected, and revised by the authors, who take full responsibility for the content of the manuscript.

\noindent G. Botelho  \\
 Instituto de Matem\'atica e Estat\'istica \\
Universidade Federal de Uberl\^andia \\
38.400-902 -- Uberl\^andia -- Brazil \\
e-mail: botelho@ufu.br

\medskip

\noindent L. A. Garcia\\
Instituto de Ci\^encias Exatas\\
Universidade Federal de Juiz de Fora\\
36.036-900 – Juiz de Fora – Brazil\\
e-mail: garcia.luis@ufjf.br

\medskip

\noindent V. C. C. Miranda\\
Departamento de Matemática\\
Instituto de Ciências Matemáticas e de Computação \\
Universidade de São Paulo \\
13566-590 -- São Carlos - SP -- Brazil  \\
e-mail: viniciusmiranda@icmc.usp.br

\end{document}